\newcount\ite\ite=1\def\0{\global\ite=1\1}
\def\1{\item{\rm(\romannumeral\the\ite)}\advance\ite1\quad}
\def\phi{\varphi}

\documentclass[12pt,twoside]{article}
\usepackage{amssymb}

\font\teneufm=eufm10 scaled \magstep1
\font\seveneufm=eufm7 scaled \magstep1
\font\fiveeufm=eufm5  scaled \magstep1
\newfam\eufmfam

\textfont\eufmfam=\teneufm
\scriptfont\eufmfam=\seveneufm
\scriptscriptfont\eufmfam=\fiveeufm

\usepackage{mathrsfs}

\newfam\msbfam
\font\tenmsb=msbm10 scaled \magstep1  \textfont\msbfam=\tenmsb
\font\sevenmsb=msbm7 scaled \magstep1 \scriptfont\msbfam=\sevenmsb
\font\fivemsb=msbm5 scaled \magstep1  \scriptscriptfont\msbfam=\fivemsb

\def\dd#1{\raise1.5pt\hbox{$\,\partial\!$}/\raise-2.5pt\hbox{$\!\partial#1\,$}}

\def\tilde{\widetilde}

\def\5#1{{\mathcal #1}}

\def\CC{{\mathbb C}}

\def\PP{{\mathbb P}}

\def\ra{\rightarrow}

\def\GL{\mathop{\rm GL}\nolimits}
\def\SL{\mathop{\rm SL}\nolimits}

\def\Ann{\mathop{\rm Ann}\nolimits}

 \def\HollowBoxx #1#2#3{{\dimen0=#1 \advance\dimen0 by -#2
       \dimen1=#1 \advance\dimen1 by #3
        \vrule height 0pt depth #3 width #2
       \hskip -#3
       \vrule height #1 depth #3 width #3}}
 \def\LeftContraction{\mathord{\kern1.45pt \HollowBoxx{6pt}{3.5pt}{.4pt}}\,}

 \def\HollowBox #1#2#3{{\dimen0=#1 \advance\dimen0 by -#3
       \dimen1=#1 \advance\dimen1 by #3
        \vrule height #1 depth #3 width #3
        \vrule height 0pt depth #3 width #2
        \hskip -#3}}
 \def\RightContraction{\mathord{\, \HollowBox{6pt}{3.1pt}{.4pt}} \kern1.6pt}

\newtheorem{theorem}{THEOREM}[section]

\newtheorem{remark}[theorem]{Remark}

\begin{document}

\begin{center}


\end{center}

\begin{center}
{\Large \bf Explicit Reconstruction of Homogeneous\\
Isolated Hypersurface Singularities\\
\vspace{0.2cm}
from their Milnor Algebras}
\medskip\medskip\\
\normalsize A. V. Isaev and N. G. Kruzhilin 
\end{center}

\begin{quotation} \small \sl\noindent  By the Mather-Yau theorem, a complex hypersurface germ ${\mathcal V}$ with isolated singularity is completely determined by its moduli algebra ${\mathcal A}({\mathcal V})$. The proof of the theorem does not provide an explicit procedure for recovering ${\mathcal V}$ from ${\mathcal A}({\mathcal V})$, and finding such a procedure is a long-standing open problem. In this paper, we present an explicit way for reconstructing ${\mathcal V}$ from ${\mathcal A}({\mathcal V})$ up to biholomorphic equivalence under the assumption that the singularity of ${\mathcal V}$ is homogeneous, in which case ${\mathcal A}({\mathcal V})$ coincides with the Milnor algebra of ${\mathcal V}$.
\end{quotation}

\thispagestyle{empty}

\pagestyle{myheadings}
\markboth{A. V. Isaev and N. G. Kruzhilin}{Reconstruction of Homogeneous Hypersurface Singularities}

\setcounter{section}{0}

\section{Introduction}\label{intro}
\setcounter{equation}{0}

Let ${\mathcal O}_n$ be the local algebra of holomorphic function germs at the origin in $\CC^n$ with $n\ge 2$. For every hypersurface germ ${\mathcal V}$ at the origin (considered with its reduced complex structure) denote by $I({\mathcal V})$ the ideal of elements of ${\mathcal O}_n$ vanishing on ${\mathcal V}$. Let $f$ be a generator of $I({\mathcal V})$, and consider the complex associative commutative algebra ${\mathcal A}({\mathcal V})$ defined as the quotient of ${\mathcal O}_n$ by the ideal generated by $f$ and all its first-order partial derivatives. The algebra ${\mathcal A}({\mathcal V})$, called the {\it moduli algebra}\, or {\it Tjurina algebra}\, of ${\mathcal V}$, is independent of the choice of $f$ as well as the coordinate system near the origin, and the moduli algebras of biholomorphically equivalent hypersurface germs are isomorphic. Clearly, ${\mathcal A}({\mathcal V})$ is trivial if and only if ${\mathcal V}$ is non-singular. Furthermore, it is well-known that $0<\dim_{\CC}{\mathcal A}({\mathcal V})<\infty$ if and only if the germ ${\mathcal V}$ has an isolated singularity (see, e.g. Chapter 1 in \cite{GLS}).

By a theorem due to Mather and Yau (see \cite{MY}), two hypersurface germs ${\mathcal V}_1$, ${\mathcal V}_2$ in $\CC^n$ with isolated singularities are biholomorphically equivalent if their moduli algebras ${\mathcal A}({\mathcal V}_1)$, ${\mathcal A}({\mathcal V}_2)$ are isomorphic. Thus, given the dimension $n$, the moduli algebra ${\mathcal A}({\mathcal V})$ determines ${\mathcal V}$ up to biholomorphism. In particular, if $\dim_{\CC}{\mathcal A}({\mathcal V})=1$, then ${\mathcal V}$ is biholomorphic to the germ of the hypersurface $\{z_1^2+\dots+z_n^2=0\}$, and if $\dim_{\CC}{\mathcal A}({\mathcal V})=2$, then ${\mathcal V}$ is biholomorphic to the germ of the hypersurface $\{z_1^2+\dots+z_{n-1}^2+z_n^3=0\}$. The proof of the Mather-Yau theorem does not provide an explicit procedure for recovering the germ ${\mathcal V}$ from the algebra ${\mathcal A}({\mathcal V})$ in general, and finding a way for reconstructing ${\mathcal V}$ (or at least some invariants of ${\mathcal V}$) from ${\mathcal A}({\mathcal V})$ is an interesting open problem (cf. \cite{Y1}, \cite{Y2}, \cite{Sch}, \cite{EI}). In this paper we present an explicit method for restoring ${\mathcal V}$ from ${\mathcal A}({\mathcal V})$ up to biholomorphic equivalence under the assumption that the singularity of ${\mathcal V}$ is homogeneous.

Let ${\mathcal V}$ be a hypersurface germ having an isolated singularity. The singularity of ${\mathcal V}$ is said to be {\it homogeneous}\, if for some (hence for every) generator $f$ of $I({\mathcal V})$ there is a coordinate system near the origin in which $f$ becomes the germ of a homogeneous polynomial. In this case $f$ lies in the Jacobian ideal ${\mathcal J}(f)$ in ${\mathcal O}_n$, which is the ideal generated by all first-order partial derivatives of $f$. Hence, for a homogeneous singularity, ${\mathcal A}({\mathcal V})$ coincides with the {\it Milnor algebra}\, ${\mathcal O}_n/{\mathcal J}(f)$ for any generator $f$ of $I({\mathcal V})$.

Let $Q(z)$, with $z:=(z_1,\dots,z_n)$, be a  holomorphic $(m+1)$-form on $\CC^n$, i.e. a homogeneous polynomial of degree $m+1$ in the variables $z_1,\dots,z_n$, where $m\ge 2$. Consider the germ ${\mathcal V}$ of the hypersurface $\{Q(z)=0\}$ and assume that: (i) the singularity of ${\mathcal V}$ is isolated, and (ii) the germ of $Q$ generates $I({\mathcal V})$. These two conditions are equivalent to the non-vanishing of the discriminant $\Delta(Q)$ of $Q$ (see Chapter 13 in \cite{GKZ}). Next, consider the gradient map ${\bf Q}:\CC^n\ra\CC^n,\,z\mapsto \hbox{grad}\,Q(z)$. Since $\Delta(Q)\ne 0$, the fiber ${\bf Q}^{-1}(0)$ consists of 0 alone; in particular, the map ${\bf Q}$ is finite at the origin. The main content of this paper is a procedure for recovering ${\bf Q}$ from ${\mathcal A}({\mathcal V})$ up to linear equivalence, where we say that two maps $\Phi_1,\Phi_2:\CC^n\ra\CC^n$ are linearly equivalent if there exist non-degenerate linear transformations $L_1,L_2$ of $\CC^n$ such that $\Phi_2=L_1\circ\Phi_1\circ L_2$. 

In fact, we consider a more general situation. Let $p_r$, $r=1,\dots,n$, be holomorphic $m$-forms on $\CC^n$ and $I$ the ideal in ${\mathcal O}_n$ generated by the germs of these forms at the origin. Define ${\mathcal A}:={\mathcal O}_n/I$ and assume that\linebreak $\dim_{\CC}{\mathcal A}<\infty$, which is equivalent to the finiteness of the map ${\bf P}:\CC^n\ra\CC^n,\,z\mapsto(p_1(z),\dots,p_n(z))$ at the origin (see Chapter 1 in \cite{GLS}). Observe that since the components of ${\bf P}$ are homogeneous polynomials, ${\bf P}$ is finite at 0 if and only if ${\bf P}^{-1}(0)=\{0\}$. In this paper we propose a procedure (which requires only linear-algebraic manipulations) for explicitly recovering the map ${\bf P}$ from ${\mathcal A}$ up to linear equivalence. As explained in Remark \ref{recognition}, this procedure also helps decide whether a given complex finite-dimensional associative algebra is isomorphic to an algebra arising from a finite homogeneous polynomial map as above.

The paper is organized as follows. Reconstruction of ${\bf P}$ from ${\mathcal A}$ is done in Section \ref{section1}.  In Section \ref{section2} we apply our method to the algebra ${\mathcal A}({\mathcal V})$ arising from $Q$ to obtain a map ${\bf Q}'$ linearly equivalent to ${\bf Q}$. It is then not hard to derive from ${\bf Q}'$ an $(m+1)$-form $Q'$ linearly equivalent to $Q$, where two forms $Q_1$, $Q_2$ on $\CC^n$ are called linearly equivalent if there exists a non-degenerate linear transformation $L$ of $\CC^n$ such that $Q_2=Q_1\circ L$. Then the germ of the hypersurface  $\{Q'(z)=0\}$ is the sought-after reconstruction of ${\mathcal V}$ from ${\mathcal A}({\mathcal V})$ up to biholomorphic equivalence. We conclude the paper by illustrating our reconstruction procedure with the example of simple elliptic singularities of type $\tilde E_6$.

{\bf Acknowledgements.} Our work was initiated during the second author's visit to the Australian National University in 2011. We gratefully acknowledge support of the Australian Research Council.

\section{Reconstruction of finite polynomial maps}\label{section1}
\setcounter{equation}{0}

Recapping the setup outlined in the introduction, let $p_r$, $r=1,\dots,n$, be holomorphic $m$-forms on $\CC^n$ and $I$ the ideal in ${\mathcal O}_n$ generated by the germs of these forms at the origin, where  $m,n\ge 2$. Define ${\mathcal A}:={\mathcal O}_n/I$ and assume that $\dim_{\CC}{\mathcal A}<\infty$ (observe that $\dim_{\CC}{\mathcal A}\ge m+1$). In this section we present a method for recovering the map ${\bf P}:\CC^n\ra\CC^n,\,z\mapsto(p_1(z),\dots,p_n(z))$ from ${\mathcal A}$ up to linear equivalence. Everywhere below we suppose that ${\mathcal A}$ is given as an abstract associative algebra, i.e. by a multiplication table with respect to some basis $e_1,\dots,e_N$, with $N:=\dim_{\CC}{\mathcal A}$. 

First of all, we find the unit ${\bf 1}$ of ${\mathcal A}$. One has ${\bf 1}=\sum_{k=1}^N\alpha_ke_k$ where the coefficients $\alpha_k\in\CC$ are uniquely determined from the linear system
$$
\sum_{k=1}^N\alpha_k(e_ke_{\ell})=e_{\ell},\quad \ell=1,\dots,N.
$$
Assume now that $e_1={\bf 1}$ and find the maximal ideal ${\mathfrak m}$ of ${\mathcal A}$. Clearly, ${\mathfrak m}$ is spanned by the vectors $e_k':=e_k-\beta_k\,{\bf 1}$, $k=2,\dots,N$, where $\beta_k\in\CC$ are uniquely fixed by the requirement that each $e_k'$ is not invertible in ${\mathcal A}$. Hence, for each $k$ the number $\beta_k$ is determined from the condition that the linear system
\begin{equation}
\sum_{\ell=1}^N\gamma_{\ell}(e_ke_{\ell}-\beta_ke_{\ell})=e_1\label{systemmax}
\end{equation}  
cannot be solved for $\gamma_1,\dots,\gamma_N\in\CC$. Since system (\ref{systemmax}) has at most one solution for any $\beta_k$, this condition is equivalent to the degeneracy of the coefficient matrix $M_k$ of (\ref{systemmax}). We have $M_k=C_k-\beta_k\hbox{Id}$, where $C_k:=(c_{kj\ell})_{j,\ell=1,\dots,N}$, with $c_{kj\ell}$ given by
$
e_ke_{\ell}=\sum_{j=1}^Nc_{kj\ell}e_j.
$
It then follows that the required value of $\beta_k$ is the (unique) eigenvalue of the matrix $C_k$.

We are now in a position to find the number of variables $n$ and the degree $m$ for the forms $p_r$. By Nakayama's lemma, ${\mathfrak m}$ is a nilpotent algebra, and we denote by $\nu$ its nil-index, which is the largest integer $\mu$ with ${\mathfrak m}^{\mu}\ne0$. Observe that $\nu\le N-1$, and therefore to determine $\nu$ it is sufficient to compute all products of the basis vectors $e_k'$ of length not exceeding $N-1$. Further, since ${\mathcal A}$ is finite-dimensional, the forms $p_r$ form a regular sequence in ${\mathcal O}_n$ (see Theorem 2.1.2 in \cite{BH}). Hence ${\mathcal A}$ is a complete intersection ring, which implies that ${\mathcal A}$ is a Gorenstein algebra (see \cite{B}). Recall that a (complex) local commutative associative algebra ${\mathcal B}$ with $1<\dim_{\CC}{\mathcal B}<\infty$ is Gorenstein if and only if for the annihilator  $\Ann({\mathfrak n}):=\{x\in{\mathfrak n}: x\,{\mathfrak n}=0\}$ of its maximal ideal ${\mathfrak n}$ one has $\dim_{\CC}\Ann({\mathfrak n})=1$ (see e.g. \cite{H}). Lemma 3.4 of \cite{Sa} yields that $\Ann({\mathfrak m})$ is spanned by the element represented by the germ of $J({\bf P}):=\det(\partial p_r/\partial z_s)_{r,s=1,\dots,n}$.

For every $i>0$, let ${\mathcal P}_i$ be the vector space of all $i$-forms on $\CC^n$ and ${\mathcal L}_i$ the linear subspace of ${\mathcal A}$ that consists of all elements represented by germs of forms in ${\mathcal P}_i$. Since ${\mathfrak m}$ consists of all elements  of ${\mathcal A}$ represented by germs in ${\mathcal O}_n$ vanishing at the origin, the subspaces ${\mathcal L}_i$ lie in ${\mathfrak m}$ and yield a grading on ${\mathfrak m}$:
$$
{\mathfrak m}=\bigoplus_{i>0} {\mathcal L}_i,\quad {\mathcal L}_i{\mathcal L}_j\subset {\mathcal L}_{i+j}\,\, \hbox{for all $i,j$.}
$$
Since $\dim_{\CC}\Ann({\mathfrak m})=1$, it immediately follows that $\Ann({\mathfrak m})={\mathfrak m}^{\nu}={\mathcal L}_d$ for $d:=\max\{i: {\mathcal L}_i\ne 0\}$. On the other hand, $\Ann({\mathfrak m})$ is spanned by the element represented by the germ of $J({\bf P})$, which is an $n(m-1)$-form. Thus $d=n(m-1)$. Furthermore, we have
\begin{equation}
\begin{array}{l}
\dim_{\CC}{\mathcal L}_i=\dim_{\CC}{\mathcal P}_i \quad\hbox{for $i=1,\dots,m-1$},\\
\vspace{-0.1cm}\\
\dim_{\CC}{\mathcal L}_m=\dim_{\CC}{\mathcal P}_m-n.
\end{array}\label{dimensions}
\end{equation} 

Now, observe that
\begin{equation}
{\mathcal L}_i={\mathcal L}_1^i\,\,\hbox{for all $i$,}\label{standard}
\end{equation}
i.e. the graded algebra ${\mathcal A}$ is {\it standard}\, in the terminology of \cite{St}. Hence ${\mathcal L}_i$ is a complement to ${\mathfrak m}^{i+1}$ in ${\mathfrak m}^i$ for all $i>0$. For $i=1$ this implies that $n$ can be recovered from the algebra ${\mathcal A}$ as 
\begin{equation}
n=\dim_{\CC}{\mathfrak m}/{\mathfrak m}^2\label{dimn}
\end{equation}
(see Remark \ref{recoverdim}). Next, for $i=\nu$ we obtain $n(m-1)=\nu$ (in particular, $n$ divides $\nu$). Thus, the degree $m$ of the forms $p_r$ can be recovered from ${\mathcal A}$ as follows: 
\begin{equation}
m=\nu/n+1.\label{formm}
\end{equation}
Note that since ${\mathcal A}$ is given as an abstract associative algebra, finding the grading $\{{\mathcal L}_i\}$ from the available data may be hard. We stress that determination of this grading is not required for recovering $n$ and $m$.

Further, choose an arbitrary basis $f_1,\dots,f_n$ in a complement to ${\mathfrak m}^2$ in ${\mathfrak m}$. Clearly, for some $C\in\GL(n,\CC)$ one has
$$
\left(
\begin{array}{c}
f_1\\
\vdots\\
f_n
\end{array}
\right)=C
\left(
\begin{array}{c}
Z_1\\
\vdots\\
Z_n
\end{array}
\right)+
\left(
\begin{array}{c}
W_1\\
\vdots\\
W_n
\end{array}
\right),
$$
where $Z_j$ is the element of ${\mathcal L}_1$ represented by the germ of the coordinate function $z_j$, and $W_j\in{\mathfrak m}^2$. Set
$$
K:=\dim_{\CC}{\mathcal P}_m=\left(
\begin{array}{c}
m+n-1\\
m
\end{array}
\right),
$$
and let $q_1(z),\dots,q_K(z)$ be the monomial basis of ${\mathcal P}_m$ where $z:=(z_1,\dots,z_n)$. Next, fix a complement ${\mathcal S}$ to ${\mathfrak m}^{m+1}$ in ${\mathfrak m}^m$ and let $\pi: {\mathfrak m}^m\ra{\mathcal S}$ be the projection onto ${\mathcal S}$ with kernel ${\mathfrak m}^{m+1}$. Condition (\ref{standard}) for $i=m$ then yields that ${\mathcal S}$ is spanned by $\pi(q_1(f)),\dots,\pi(q_K(f))$, where $f:=(f_1,\dots,f_n)$. On the other hand, by (\ref{dimensions}) we have $\dim_{\CC}{\mathcal S}=K-n$. Hence one can find $n$ linear relations
\begin{equation}
\sum_{\rho=1}^K\gamma_{\sigma\rho}\pi(q_{\rho}(f))=0,\quad \sigma=1,\dots,n,\label{relations}
\end{equation}
where the vectors $\gamma_{\sigma}:=(\gamma_{\sigma\,1},\dots,\gamma_{\sigma\,K})\in\CC^K$ are linearly independent. 

Further, extracting from (\ref{relations}) the ${\mathcal L}_m$-components one obtains
\begin{equation}
\sum_{\rho=1}^K\gamma_{\sigma\rho}q_{\rho}(CZ)=0,\quad \sigma=1,\dots,n,\label{relations1}
\end{equation}
where $Z:=(Z_1,\dots,Z_n)$. Identity (\ref{relations1}) is equivalent to 
\begin{equation}
\sum_{\rho=1}^K\gamma_{\sigma\rho}q_{\rho}(Cz)\in I,\quad \sigma=1,\dots,n,\label{relations2}
\end{equation}
where each $q_{\rho}(z)$ is identified with its germ at the origin. From (\ref{relations2}) one immediately obtains that for some matrix $D\in\GL(n,\CC)$ the following holds:
$$
\Gamma q(Cz)\equiv D{\bf P}(z),
$$
where $\Gamma:=(\gamma_{\sigma\rho})_{\sigma=1,\dots,n,\,\rho=1,\dots,K,}$ and $q:=(q_1,\dots,q_K)$. Thus, the map
\begin{equation}
\Phi: \CC^n\ra\CC^n,\quad z\mapsto \Gamma q(z)\label{recmap}
\end{equation}
is linearly equivalent to ${\bf P}$ as required.

We now summarize the main steps of our algorithm for recovering ${\bf P}$ from ${\mathcal A}$ up to linear equivalence:
$$
\begin{array}{ll}
1.& \hbox{Find ${\mathfrak m}$ and its nil-index $\nu$.}\\
\vspace{-0.4cm}\\
2.& \hbox{Determine $n$ from formula (\ref{dimn}).}\\
\vspace{-0.4cm}\\
3.& \hbox{Determine $m$ from formula (\ref{formm}).}\\
\vspace{-0.4cm}\\
4.& \hbox{Choose a complement to ${\mathfrak m}^2$ in ${\mathfrak m}$ and an arbitrary basis $f_1,\dots,f_n$ in}\\
&\hbox{this complement.}\\
\vspace{-0.4cm}\\
5. & \hbox{Calculate $q_1(f),\dots,q_K(f)$, where $f:=(f_1,\dots,f_n)$ and $q_1(z),\dots,q_K(z)$}\\
&\hbox{are all monomials of degree $m$ in $z:=(z_1,\dots,z_n)$.}\\
\vspace{-0.4cm}\\
6. & \hbox{Choose a complement ${\mathcal S}$ to ${\mathfrak m}^{m+1}$ in ${\mathfrak m}^m$.}\\
\vspace{-0.4cm}\\
7. & \hbox{Compute $\pi(q_1(f)),\dots,\pi(q_K(f))$, where $\pi: {\mathfrak m}^m\ra{\mathcal S}$ is the projection}\\
&\hbox{onto ${\mathcal S}$ with kernel ${\mathfrak m}^{m+1}$.}\\
\vspace{-0.4cm}\\
8. & \hbox{Find $n$ linearly independent linear relations among the vectors}\\
&\hbox{$\pi(q_1(f)),\dots,\pi(q_K(f))$ as in (\ref{relations}).}\\
\vspace{-0.4cm}\\
9. & \hbox{Formula (\ref{recmap}) then gives a map linearly equivalent to ${\bf P}$.}
\end{array}
$$
In the next section this algorithm will be applied to the gradient map arising from a form with non-zero discriminant.

\begin{remark}\label{recognition} \rm A natural problem is to characterize the algebras that arise from finite homogeneous polynomial maps as above among all complex finite-dimensional Gorenstein algebras. This problem is a special case of the well-known {\it recognition problem}\, for the moduli algebras of general isolated hypersurface singularities and the corresponding Lie algebras of derivations (see, e.g. \cite{Y1}, \cite{Y2}, \cite{Sch}). The algorithm presented here can help decide whether a given finite-dimensional Gorenstein algebra ${\mathcal B}$ is isomorphic to an algebra ${\mathcal A}$ of the kind considered in this section. Indeed, one can attempt to formally apply the algorithm to ${\mathcal B}$. For the algorithm to go through one requires that: (i) the nil-index of the maximal ideal of ${\mathcal B}$ be divisible by the number $n$ found from formula (\ref{dimn}), (ii) for some basis $f_1,\dots,f_n$ in some complement to ${\mathfrak m}^2$ in ${\mathfrak m}$ and for some complement ${\mathcal S}$ to ${\mathfrak m}^{m+1}$ in ${\mathfrak m}^m$ there exist $n$ linearly independent linear relations among the vectors $\pi(q_1(f)),\dots,\pi(q_K(f))$, with $m$ being the number found from formula (\ref{formm}), and (iii) the map $\Phi:\CC^n\ra\CC^n$ produced on Step 9 be finite at the origin. If the algorithm fails (i.e. some of conditions (i)--(iii) are not satisfied), then ${\mathcal B}$ does not arise from a finite homogeneous polynomial map. If the algorithm successfully finishes, the resulting map $\Phi$ is a candidate map from which ${\mathcal B}$ may potentially arise. In order to see if this is indeed the case, one needs to check whether ${\mathcal B}$ is isomorphic to the algebra associated to $\Phi$. For this purpose one can use the criterion for isomorphism of finite-dimensional Gorenstein algebras established in \cite{FIKK}.          
\end{remark}

\section{Reconstruction of homogeneous\\ singularities}\label{section2}
\setcounter{equation}{0}

Suppose now that ${\bf P}={\bf Q}:=\hbox{grad}\,Q$ for a holomorphic $(m+1)$-form $Q$ on $\CC^n$ with $\Delta(Q)\ne 0$, where $\Delta$ is the discriminant. Let $\Phi$ be a map linearly equivalent to ${\bf Q}$ produced by the procedure described in Section \ref{section1} from the algebra ${\mathcal A}({\mathcal V})$, where ${\mathcal V}$ is the germ of the hypersurface $\{Q=0\}$ at the origin. We then have
\begin{equation}
\Phi(z)\equiv C_1\,\hbox{grad}\,Q(C_2z)\label{solut1}
\end{equation}
for some $C_1,C_2\in\GL(n,\CC)$. Our next task is to recover $Q$ from $\Phi$ up to linear equivalence. 

Let $Q'$ be the $(m+1)$-form defined by $Q'(z):=Q(C_2z)$ for all $z\in\CC^n$. Then $\hbox{grad}\,Q(C_2z)=(C_2^{-1})^T\hbox{grad}\,Q'(z)$, and (\ref{solut1}) implies
$$   
\Phi(z)\equiv C\,\hbox{grad}\,Q'(z)
$$
for some $C\in\GL(n,\CC)$. For any $n\times n$-matrix $D$ we now introduce the holomorphic differential 1-form $\omega^D:=\sum_{r=1}^n\Psi^D_rdz_r$ on $\CC^n$, where $(\Psi^D_1,\dots,\Psi^D_n)$ are the components of the map $\Psi^D:=D\,\Phi$. Consider the equation
\begin{equation}
d\omega^D=0\label{system}
\end{equation}
as a linear system with respect to the entries of the matrix $D$. Clearly, $C^{-1}$ is a solution of (\ref{system}). Let $D_0$ be another solution of (\ref{system}) and assume that $D_0\in\GL(n,\CC)$. Every closed holomorphic differential form on $\CC^n$ is exact, and integrating $\Psi^{D_0}$ one obtains an $(m+1)$-form $Q''$ on $\CC^n$. Then $\hbox{grad}\,Q''=\Psi^{D_0}=D_0C\,\hbox{grad}\,Q'$, and therefore $\Delta(Q'')\ne 0$. Furthermore, the Milnor algebras of the germs ${\mathcal V}'$ and ${\mathcal V}''$ of the hypersurfaces $\{Q'(z)=0\}$ and $\{Q''(z)=0\}$ coincide. By the Mather-Yau theorem, this implies that ${\mathcal V}'$ and ${\mathcal V}''$ are biholomorphically equivalent and therefore ${\mathcal V}''$ is biholomorphically equivalent to ${\mathcal V}$, which yields that $Q''$ is linearly equivalent to $Q$. Thus, any non-degenerate matrix that solves linear system (\ref{system}) leads to an $(m+1)$-form linearly equivalent to $Q$ and a hypersurface germ biholomorphically equivalent to ${\mathcal V}$.

We will now illustrate our method for recovering ${\mathcal V}$ from ${\mathcal A}({\mathcal V})$ by the example of simple elliptic singularities of type $\tilde E_6$. These singularities form a family parametrized by $t\in\CC$ satisfying $t^3+27\ne 0$. Namely, for every such $t$ let ${\mathcal V}_t$ be the germ at the origin of the hypersurface $\{Q_t(z)=0\}$, where $Q_t$ is the following cubic on $\CC^3$:
$$
Q_t(z):=z_1^3+z_2^3+z_3^3+tz_1z_2z_3,\quad z:=(z_1,z_2,z_3).
$$
Below we explicitly show how $Q_t$ can be recovered from the algebra\linebreak ${\mathcal A}_t:={\mathcal A}({\mathcal V}_t)$ up to linear equivalence. 

Recall that the starting point of our reconstruction procedure is a multiplication table with respect to some basis. The algebra ${\mathcal A}_t$ has dimension 8 and with respect to a certain basis $e_1,\dots,e_8$ is given by (see Remark \ref{basis} below):
\begin{equation}
\begin{array}{l}
e_1e_j=e_j\,\,\hbox{for $j=1,\dots,8$},\,e_2^2=\displaystyle-\frac{t}{3}e_3+\frac{2t}{3}e_6,\, e_2e_3=e_6,\\
\vspace{-0.1cm}\\
e_2e_4=e_5-e_6-e_8,\, e_2e_5=e_7,\,e_2e_6=0,\, e_2e_7=0,\,e_2e_8=e_7,\\
\vspace{-0.1cm}\\
e_3e_j=0\,\,\hbox{for $j=3,\dots,8$},\,e_4^2=\displaystyle -\frac{t}{3}e_7,\,e_4e_5=e_3-2e_6,\,e_4e_6=0,\\
\vspace{-0.1cm}\\
e_4e_7=e_6,\,e_4e_8=e_3-2e_6,\,e_5^2=\displaystyle -\frac{t}{3}e_5+(2+t)e_6+\frac{t}{3}e_8,\,e_5e_6=0,\\
\vspace{-0.1cm}\\
e_5e_7=0,\,e_5e_8=\displaystyle -\frac{t}{3}e_5+(1+t)e_6+\frac{t}{3}e_8,\,e_6e_j=0\,\,\hbox{for $j=6,7,8$},\\
\vspace{-0.1cm}\\
e_7e_j=0\,\,\hbox{for $j=7,8$},\,e_8^2=\displaystyle-\frac{t}{3}e_5+te_6+\frac{t}{3}e_8.
\end{array}\label{table}
\end{equation}
It is clear from (\ref{table}) that $e_1={\bf 1}$ and ${\mathfrak m}_t=\langle e_2,\dots,e_8\rangle$, where $\langle {}\cdot{} \rangle$ denotes linear span and ${\mathfrak m}_t$ is the maximal ideal of ${\mathcal A}_t$. We then have ${\mathfrak m}_t^2=\langle e_3, e_6, e_7,e_5-e_8\rangle$, ${\mathfrak m}_t^3=\langle e_6\rangle$, ${\mathfrak m}_t^4=0$, hence $\nu=3$. Further, by formula (\ref{dimn}) we obtain $n=3$, which together with formula (\ref{formm}) yields $m=2$.

We now list all monomials of degree 2 in $z$ as follows:
$$
q_1(z):=z_1^2,\,q_2(z):=z_2^2,\,q_3(z):=z_3^2,\,q_4(z):=z_1z_2,\,q_5(z):=z_1z_3,\,q_6(z):=z_2z_3
$$
(here $K=6$). Next, we let $f_1:=e_2$, $f_2:=e_4$, $f_3:=e_5$, which for $f:=(f_1,f_2,f_3)$ yields
$$
\begin{array}{l}
q_1(f)=\displaystyle-\frac{t}{3}e_3+\frac{2t}{3}e_6,\,q_2(f)=-\frac{t}{3}e_7,\,q_3(f)=\displaystyle -\frac{t}{3}e_5+(2+t)e_6+\frac{t}{3}e_8,\\
\vspace{-0.1cm}\\
q_4(f)=e_5-e_6-e_8,\quad q_5(f)=e_7,\quad q_6(f)=e_3-2e_6.
\end{array}
$$

Further, define ${\mathcal S}:=\langle e_3,e_7,e_5-e_8\rangle$. Clearly, ${\mathcal S}$ is a complement to ${\mathfrak m}_t^3$ in ${\mathfrak m}_t^2$. Then for the projection $\pi: {\mathfrak m}_t^2\ra{\mathcal S}$ with kernel ${\mathfrak m}_t^3$ one has
$$
\begin{array}{l}
\pi(q_1(f))=\displaystyle-\frac{t}{3}e_3,\,\pi(q_2(f))=-\frac{t}{3}e_7,\,\pi(q_3(f))=-\frac{t}{3}e_5+\frac{t}{3}e_8,\\
\vspace{-0.1cm}\\
\pi(q_4(f))=e_5-e_8,\quad \pi(q_5(f))=e_7,\quad \pi(q_6(f))=e_3.
\end{array}
$$
The vectors $\pi(q_1(f)),\dots,\pi(q_6(f))$ satisfy the following three linearly independent linear relations:
$$
\pi(q_1(f))+\frac{t}{3}\pi(q_6(f))=0,\,\,\pi(q_2(f))+\frac{t}{3}\pi(q_5(f))=0,\,\,\pi(q_3(f))+\frac{t}{3}\pi(q_4(f))=0.
$$
Hence we have 
$$
\Gamma=\left(
\begin{array}{cccccc}
1 & 0 & 0 & 0 & 0 & t/3\\
0 & 1 & 0 & 0 & t/3 & 0\\
0 & 0 & 1 & t/3 & 0 & 0
\end{array}
\right),
$$
which for $q(z):=(q_1(z),\dots,q_6(z))$ yields
$$
\Phi(z)=\Gamma q(z)=\left(
\begin{array}{l}
z_1^2+\displaystyle\frac{t}{3}z_2z_3\\
\vspace{-0.1cm}\\
z_2^2+\displaystyle\frac{t}{3}z_1z_3\\
\vspace{-0.1cm}\\
z_3^2+\displaystyle\frac{t}{3}z_1z_2\\
\vspace{-0.1cm}\\
\end{array}
\right).
$$

It remains to recover $Q_t$ from $\Phi$ up to linear equivalence. For 
$$
D=\left(
\begin{array}{lll}
d_{11} & d_{12} & d_{13}\\
d_{21} & d_{22} & d_{23}\\
d_{31} & d_{32} & d_{33}
\end{array}
\right)
$$
system (\ref{system}) is equivalent to the following system of equations:
\begin{equation}
\begin{array}{lll}
\displaystyle 2d_{12}-\frac{t}{3}d_{23}=0,&\displaystyle 2d_{21}-\frac{t}{3}d_{13}=0, & td_{11}-td_{22}=0,\\
\vspace{-0.1cm}\\
\displaystyle 2d_{13}-\frac{t}{3}d_{32}=0,&\displaystyle 2d_{31}-\frac{t}{3}d_{12}=0, & td_{11}-td_{33}=0,\\
\vspace{-0.1cm}\\
\displaystyle 2d_{23}-\frac{t}{3}d_{31}=0,&\displaystyle 2d_{32}-\frac{t}{3}d_{21}=0, & td_{22}-td_{33}=0.  
\end{array}\label{sys1}
\end{equation}
If $t\ne 0$ and $t^3\ne 216$, the only non-degenerate solutions of (\ref{sys1}) are non-zero scalar matrices. Integrating $\Psi^D$ for such a matrix $D$ we obtain a form proportional to $Q_t$, which is obviously linearly equivalent to $Q_t$. If $t=0$, any non-degenerate solution of (\ref{sys1}) is a diagonal matrix with non-zero $d_{11}$, $d_{22}$, $d_{33}$. Integrating $\Psi^D$ for such a matrix $D$ we obtain the form
$$
\frac{1}{3}\left(d_{11}z_1^3+d_{22}z_2^3+d_{33}z_3^3\right),
$$
which is linearly equivalent to $Q_0=z_1^3+z_2^3+z_3^3$ by suitable dilations of the variables.

The remaining case $t^3=216$ is more interesting. Writing $t=6\lambda$ with $\lambda^3=1$, we see that $D$ is a solution of (\ref{sys1}) if and only if
$$
d_{11}=d_{22}=d_{33},\, d_{12}=\lambda^2d_{31},\, d_{23}=\lambda d_{31},\,d_{21}=\lambda^2 d_{32},\,d_{13}=\lambda d_{32}.
$$
Such a matrix $D$ is non-degenerate if and only if $d_{11}^3+d_{31}^3+d_{32}^3-3\lambda d_{11}d_{31}d_{32}\ne 0$. For example, letting $d_{11}=0$, $d_{31}=0$, $d_{32}=1$ one obtains
$$
\Psi^D=\left(
\begin{array}{l}
\lambda z_3^2+2\lambda^2 z_1z_2\\
\vspace{-0.1cm}\\
\lambda^2 z_1^2+2 z_2z_3\\
\vspace{-0.1cm}\\
z_2^2+2\lambda z_1z_3
\end{array}
\right).
$$    
Integration of $\Psi^D$ leads to the form ${\mathcal Q}_{\lambda}:=\lambda^2 z_1^2z_2+\lambda z_1 z_3^2+z_2^2z_3$. As we have noted above, the Mather-Yau theorem implies that ${\mathcal Q}_{\lambda}$ is linearly equivalent to $Q_t$. Furthermore, the cubic ${\mathcal Q}_{\lambda}$ is equivalent to ${\mathcal Q}_1$ by the map $(z_1,z_2,z_3)\mapsto(z_1/\lambda,z_2,z_3)$. Hence each of the three cubics ${\mathcal Q}_{\lambda}$ with $\lambda^3=1$ is linearly equivalent to each of the three cubics $Q_t$ with $t^3=216$. 

This last fact can also be understood without referring to the Mather-Yau theorem, as follows. It is well-known that all non-equivalent ternary cubics with non-vanishing discriminant are distinguished by the invariant   
$$
{\tt J}:=\frac{{\tt I}_4^3}{\Delta},
$$
where ${\tt I}_4$ is a certain classical $\SL(3,\CC)$-invariant of degree 4 (see, e.g. pp. 381--389 in \cite{El}). For any ternary cubic $Q$ with $\Delta(Q)\ne 0$ one has ${\tt J}(Q)=j(Z_Q)/110592$ where $j(Z_Q)$ is the value of the $j$-invariant for the elliptic curve $Z_Q$ in $\CC\PP^2$ defined by $Q$. Details on computing ${\tt J}(Q)$ for any $Q$ can be found, for example, in \cite{Ea}. In particular, ${\tt J}({\mathcal Q}_1)=0$ and for the cubic $Q_t$ with any $t\in\CC$, $t^3+27\ne 0$, one has
$$
{\tt J}(Q_t)=-\frac{t^3(t^3-216)^3}{110592(t^3+27)^3}.
$$
It then follows that each of the cubics ${\mathcal Q}_{\lambda}$ with $\lambda^3=1$ is linearly equivalent to each of the cubics $Q_t$ with $t^3=216$, as stated above.  

\begin{remark}\label{basis} \rm One basis in which the algebra ${\mathcal A}_t$ is given by multiplication table (\ref{table}) is as follows:
$$
\begin{array}{l}
e_1={\bf 1},\,e_2=Z_1+Z_1Z_3,\,e_3=Z_2Z_3+2Z_1Z_2Z_3,\,e_4=Z_2+Z_2Z_3,\\
\vspace{-0.1cm}\\
e_5=Z_3+Z_1Z_2+3Z_1Z_2Z_3,\,e_6=Z_1Z_2Z_3,\,e_7=Z_1Z_3,\,e_8=Z_3.
\end{array}
$$
Note that in our reconstruction of $Q_t$ from ${\mathcal A}_t$ above we only used table (\ref{table}), not the explicit form of the basis.
\end{remark}

\begin{remark}\label{recoverdim} \rm As was noted in the introduction, a hypersurface germ ${\mathcal V}$ with isolated singularity is determined, in general, by the algebra ${\mathcal A}({\mathcal V})$ and the dimension $n$ of the ambient space. We stress that in the case of homogeneous singularities the dimension $n$ can be extracted from ${\mathcal A}({\mathcal V})$ (see formula (\ref{dimn})).    
\end{remark}

{\obeylines
Department of Mathematics
The Australian National University
Canberra, ACT 0200
Australia
e-mail: alexander.isaev@anu.edu.au
\hbox{ \ \ }
Department of Complex Analysis
Steklov Mathematical Institute
8 Gubkina St.
Moscow GSP-1 119991
Russia
e-mail: kruzhil@mi.ras.ru
}


\begin{thebibliography}{ABCD}

\bibitem[B]{B} Bass, H., On the ubiquity of Gorenstein rings, {\it Math. Z.} 82 (1963), 8--28.

\bibitem[BH]{BH} Bruns, W. and Herzog, J., {\it Cohen-Macaulay Rings}, Cambridge Studies in Advanced Mathematics, 39, Cambridge University Press, Cambridge, 1993.


\bibitem[Ea]{Ea} Eastwood, M. G., Moduli of isolated hypersurface singularities, {\it Asian J. Math.} 8 (2004), 305--313.

\bibitem[EI]{EI} Eastwood, M. G. and Isaev, A. V., Extracting invariants of isolated hypersurface singularities from their moduli algebras, preprint, available from http://arxiv.org/abs/1110.2559.

\bibitem[El]{El} Elliott, E. B., {\it An Introduction to the Algebra of Quantics}, Oxford University Press, 1895.

\bibitem[FIKK]{FIKK} Fels, G., Isaev, A., Kaup, W. and Kruzhilin, N., Isolated hypersurface singularities and special polynomial realizations of affine quadrics, {\it J. Geom. Analysis} 21 (2011), 767--782.

\bibitem[GKZ]{GKZ} Gelfand, I. M., Kapranov, M. M. and Zelevinsky, A. V., {\it Discriminants, Resultants and Multidimensional Determinants}, Modern Birkh\"auser Classics, Birkh\a"user Boston, Inc., Boston, MA, 2008.

\bibitem[GLS]{GLS} Greuel, G.-M., Lossen, C. and Shustin, E., {\it Introduction to Singularities and Deformations}, Springer Monographs in Mathematics, Springer, Berlin, 2007.

\bibitem[H]{H} Huneke, C., Hyman Bass and ubiquity: Gorenstein rings, in {\it Algebra, $K$-theory, Groups, and Education} (New York, 1997), Contemp. Math., 243, Amer. Math. Soc., Providence, RI, 1999, pp. 55--78.


\bibitem[MY]{MY} Mather, J. and Yau, S. S.-T., Classification of isolated hypersurface singularities by their moduli algebras, {\it Invent. Math.} 69 (1982), 243--251.

\bibitem[Sa]{Sa} Saito, K., Einfach-elliptische Singularit\"aten, {\it Invent. Math.} 23 (1974), 289--325.

\bibitem[Sch]{Sch} Schulze, M., A solvability criterion for the Lie algebra of derivations of a fat point, {\it J. Algebra} 323 (2010), 2916--2921.

\bibitem[St]{St} Stanley, R., Hilbert functions of graded algebras, {\it Advances in Math.} 28 (1978), 57--83.

\bibitem[Y1]{Y1} Yau, S. S.-T., Solvable Lie algebras and generalized Cartan matrices arising from isolated singularities, {\it Math. Z.} 191 (1986), 489--506.

\bibitem[Y2]{Y2} Yau, S. S.-T., Solvability of Lie algebras arising from isolated singularities and nonisolatedness of singularities defined by ${\mathfrak {sl}}(2,\CC)$ invariant polynomials, {\it Amer. J. Math.} 113 (1991), 773--778.

\end{thebibliography}
\end{document}